\documentclass[a4paper, 11pt]{article}

\usepackage{graphicx}
\usepackage{amsmath}
\usepackage{amsfonts}
\usepackage{amssymb}
\usepackage{authblk}

\usepackage{tikz}
\usetikzlibrary{arrows,backgrounds,snakes,patterns}
\usetikzlibrary{shapes,arrows,chains}
\usepackage{verbatim}
\usepackage{booktabs}

\usepackage[flushleft]{threeparttable}

\begin{document}

\title{A Decomposition Approach to Solve The Quay Crane Scheduling Problem}

\author[1]{Afonso Sampaio \footnote{afonsohs@dcc.ufmg.br}}
\author[2]{Sebasti\'{a}n Urrutia \footnote{surrutia@dcc.ufmg.br}}
\author[1]{Johan Oppen \footnote{johan.oppen@himolde.no}}

\affil[1]{Molde University College, Molde, Norway}
\affil[2]{Universidade Federal de Minas Gerais, Belo Horizonte, Brazil}

\maketitle

\begin{abstract}
%In this work, we propose a decomposition approach to exactly solve the Quay 
In this work we propose a decomposition approach to solve the quay 
crane scheduling problem. This is an important maritime transportation 
problem faced in container terminals where quay cranes are used to handle 
cargo. The objective is to determine a sequence of loading and unloading 
operations for each crane in order to minimize the completion time. We 
solve a mixed integer programming formulation for the quay 
crane scheduling problem, decomposing it into a vehicle routing 
problem and a corresponding scheduling problem. The routing sub-problem is 
solved by minimizing the longest crane completion time without taking crane 
interference into account.
This solution provides a lower bound for the 
makespan of the whole problem and is sent to the scheduling sub-problem, where 
a completion time for each task and the makespan are determined. 
%This information is used to update the best solution found and to provide 
%cuts to 
%the routing problem. Cuts are used to avoid the generation of the same or 
%similar routes that cannot further improve the best known solution.
This scheme resembles Benders' decomposition and, in particular, the scheme 
underlying combinatorial Benders' cuts.
%, but the cut generation we propose does 
%not rely on finding an irreducible infeasible subsystem. Rather, the 
%scheduling 
%sub-problem is always feasible and we derive cuts from the scheduling 
%solution. 
We evaluate the proposed approach by solving instances from the literature and comparing the 
results with other available methods.

%since the time required to load and unload cargo depends on the cranes assigned to operate in the vessel.
%Side constraints such as precedence and non-simultaneity relationships among tasks and interference situations of cranes must be taken into account.
% considering the processing times of the tasks and the time to move from one bay to another
%The solution, which provides a lower bound for the vessel makespan, determines the assignment and the ordering of tasks for each crane and is sent to the scheduling sub-problem, where a completion time for each task and the makespan are determined.

%\keywords{Container terminal \and Crane scheduling \and Decomposition methods}
\end{abstract}

\section{Introduction}
\label{sec:intro}
%In a sea port, maritime container terminals provide logistic facilities for the transshipment of cargo between container ships (\review{vessels=ships?}) and land vehicles. In the quayside area, loading and unloading of the vessels are performed by quay cranes and technical performance of a terminal is often measured in terms of how the usage of these equipments maximize the time for completing one ship's operations or minimize all ship's delays in the port.
 
Container terminals in ports provide logistic facilities for transshipment of 
cargo between container vessels and other modes of transportation. In the 
quayside area, loading and unloading of vessels are performed by quay cranes. 
The performance of a terminal is often measured in terms of how efficiently 
these cranes are utilized to minimize loading and unloading times or to minimize 
the delay of ships in the port. 

%Container terminals operations for shipping planning are (\review{often?}) comprised by the berth planning, the quay crane assignment and the quay crane scheduling problem (QCSP). The berth in which a vessel moors at the port has to be allocated before the ship's arrival. The number and which kind of quay cranes assigned to a vessel depend on the particular berth allocated, the technical data of ships and quay cranes, and the contracts between the terminal and shipping companies. These three problems are highly integrated and the optimal coordination of all these processes is very complex, specially for large container terminals. To deal with this situation, each problem is tackled separately at the cost of reducing the quality of the solution \review{(improve this...)}.

Quay side operations are comprised by the berth allocation, the quay crane assignment and the quay crane scheduling problem (QCSP). The berth in which a vessel moors at the port has to be allocated before the ship's arrival. Number and types of quay cranes assigned to a vessel depend on the particular berth allocated, the technical data of ships and quay cranes, and the contracts between the terminal and shipping companies. These three problems are highly integrated and should ideally be handled together. Most of the research literature does, however, tackle these problems separately, as combining them would lead to a very complex planning problem.

In this paper, we tackle the QCSP. Given the berth and cranes allocated to serve 
a given container vessel, the QCSP aims at determining a schedule for the cranes 
to handle the unloading and loading tasks on the sections (bays) of the vessel. 
Each task (loading or unloading of a container or a group of containers) must be 
assigned to a crane, and the sequence of tasks performed by the cranes has to be 
determined such as to minimize the total time the ship spends at the terminal. 
Since the cranes are mounted on the same rail, crossing situations must be 
avoided in the schedule. In addition, adjacent cranes have to keep a safety 
distance between them at all times.

In the QCSP, the tasks to be scheduled can be described according to the level 
of granularity considered in the model. A task consists of the loading or 
unloading of a single or a group of containers in certain areas of the vessel. 
In this paper, we consider the QCSP with container groups, in which a task 
consists of the loading/unloading of a group of containers in a bay, and 
container groups in the same bay can be assigned to different cranes. For a more 
thorough description of other classifications, the reader is referred to 
\cite{bier2010}. Precedence relationships between two tasks are derived from 
their relative positions in the ship and the operation type. Within the same 
bay, unloading always precedes loading, an unloading operation from the deck 
must be performed before an unloading operation from the hold, and loading in 
the hold must be performed before loading on the deck. The non-simultaneity of 
tasks is the result of the safety distance that must be kept between two cranes, 
and as a way to avoid workload peaks in some areas 
in the 
terminal.

The rest of the paper is organized as follows. In the next section we survey some of the relevant literature on the QCSP. In Section \ref{sec:math}, a mixed integer programming (MIP) formulation for the problem is presented together with the revisions proposed to a correct treatment of crane interference. Our decomposition algorithm is presented in section \ref{sec:alg} and we compare it with a MIP model and other algorithmic approaches in Section \ref{sec:resul}.

\section{Literature review}
\label{sec:lit}
For an extensive review on container terminal operations, the reader is referred to \cite{stah2008} and \cite{dirk2004}. A survey on berth allocation and quay crane scheduling problems is provided by \cite{bier2010}. 

The QCSP with container groups was first addressed in \cite{daganzo1989}. In 
that work, a MIP formulation for the problem was 
presented. The objective was to minimize the weighted sum of departure times of 
vessels. Container groups belonging to the same bay could be assigned to 
different cranes, but the model did not take crane interference into 
consideration. In \cite{peterkofsky1990}, a branch-and-bound is developed to 
tackle larger instances than those addressed in \cite{daganzo1989}.

More recently, another MIP formulation using big-M values for the QCSP with container groups was presented by \cite{kim2004}, in which cranes serve a single vessel. Crane operations take into account precedence relations between load and unload tasks and interference among cranes, but only instances involving up to two quay cranes and six tasks could be solved to optimality. A branch-and-bound method was proposed  along with a local search GRASP heuristic to tackle larger instances. 

An improved MIP model was proposed by \cite{moccia2006}, who noticed that the 
model of \cite{kim2004} does not correctly address cranes interference. Besides 
providing a revised model, the authors proposed lower bounds for 
task starting times and an upper bound for the task completion time, and used 
these values to improve the big-M values used in the formulation. A 
branch-and-cut method using new valid inequalities and other inequalities 
adapted from the Precedence Constrained Travelling Salesman Problem \cite{balas1995} was presented, and 
significant improvements were found for the benchmark instances from 
\cite{kim2004}.

In \cite{samarra2007}, the problem is decomposed into a routing and a scheduling 
problem. The authors use a Tabu Search heuristic for the routing problem, 
whereas a local search using a neighborhood defined over a disjunctive graph is 
used to handle the scheduling part. The algorithm outperforms the GRASP 
algorithm in \cite{kim2004} and produces slightly weaker solutions at the 
expense of better running times when compared to the branch-and-cut presented 
in \cite{moccia2006}.

The authors of \cite{bier2009} noticed situations where crossing between the 
cranes is not detected by the previous MIP models. To correct these situations, 
they introduced a suitable temporal distance between any two tasks performed by 
two different cranes so that this time period allows the tasks to be performed 
without cranes crossing. The authors propose a heuristic algorithm based on a 
branch-and-bound method searching in a reduced solution space consisting of 
unidirectional schedules (i.e. schedules in which cranes move in one 
direction only, either left to right or right to left). Considering the set of benchmark instances with two and three cranes, new best solutions were found and 
improved execution times were obtained.

%\review{
The idea of exclusively search the space of unidirectional schedules is further explored in \cite{legato2012} and \cite{chen2014}. In the former work, the authors enrich the traditional QCSP model with other aspects of practical relevance such as cranes with non-uniform productivity rates, time windows and independent unidirectionality (a schedule in which each crane moves along the same direction when serving the vessel, but the directions can be different for each crane). Based on the branch-and-bound method in \cite{bier2009}, the authors propose an approach to take into account these new features, except for the so-called independent unidirectional schedule. Benchmark instances from the literature and instances from the port of Gioia Tauro, Italy are solved. In the latter work, the authors propose a novel mathematical model for the unidirectional QCSP with container groups (cluster-based). Whereas most of the MIP formulations for the QCSP consider decision variables for the sequence of assigned tasks to each crane and for the relative order in which tasks are processed, in most of the cases the number of tasks is relatively larger when compared to the number of cranes and bays. The proposed model extends the MIP model in \cite{Liu2006}, considering binary decision variables only for the assignment of tasks to cranes. By introducing fewer binary variables and considering only unidirectional schedules, the authors obtain an easy-to-formulate model which can be quickly solved by off-the-shelf optimizers. The results show that this approach is superior to the algorithm proposed by \cite{legato2012}. Also, the authors have identified a possible blocking phenomena at the beginning of the scheduling when the ready times are non-zero and proposed a way to correct the previous formulations.

As can be seen in \cite{bier2009}, \cite{legato2012} and \cite{chen2014}, the exclusive consideration of unidirectional schedules is an effective strategy for obtaining good solutions for the QCSP. Nevertheless, such schedules may be of limited use in practical situations: stability might be an issue if all cranes operate heavy material at the same time on one end of the vessel \cite{bier2009}. Moreover, research towards identifying the existence of optimality conditions for unidirectional schedules can be improved with exact approaches for the QCSP without restricting the search space.
%}

\section{Mathematical models}
\label{sec:math}
The MIP formulation for the QCSP used in this work is based on the model 
provided by \cite{kim2004}, the modifications  and enhancements reported by 
\cite{moccia2006} and on the revision introduced in \cite{bier2009} to include a 
correct treatment of crane interference. The problem is defined over a set of 
tasks $\Omega=\{1,...,n\}$ and a set of cranes $K=\{1,...,q\}$. Let $0$ and $T$ 
be artificial tasks of null processing time that model the crane at its initial 
state (i.e. before it has executed any task) and at 
its final state (i.e. after it has completed its assigned tasks), 
respectively. Let $\Omega^0=\Omega \cup \{0\}$ and $\Omega^T=\Omega \cup \{T\}$.
The vessel is divided in $B$ bays, each task $i\in \Omega$ has a 
processing time $p_i$ and a bay number $l_i\in 
\mathbb{Z}_+,~l_i\leq B$, representing where the task is 
located. The inter-crane safety margin is modelled by $\delta$, and between any two cranes there must be $\delta$ bays without cranes at any time. Precedence relationships between tasks are expressed 
by the set $\Phi=\{(i,j)|i,j\in \Omega, i\prec j\}$, where $i\prec j$ means task 
$i$ must be completed before task $j$ starts. Non simultaneity of tasks are 
handled by the set $\Psi=\{(i,j)|i,j\in \Omega, i\nparallel j\}$, where $i 
\nparallel j$ means that tasks $i$ and $j$ cannot be performed concomitantly. 
Note that if tasks $i$ and $j$ are in the same bay (i.e. $l_i=l_j$) then 
$(i,j)\in \Psi$. Also, we have $\Phi \subseteq \Psi$

Each crane $k$ has an earliest available time, $r_k$, and it is initially 
positioned at bay $l_0^k$. Also, after the crane has completed all the assigned 
tasks, it must be positioned at bay $l_T^k$. If $l_T^k=0$, then the final position of the crane does not 
matter. Let $t$ be the time of travelling 
between two adjacent bays. The time needed for a crane to move from the bay 
where task $i$ is located to the bay where task $j$ is located is 
$t_{ij}=t\times |l_i-l_j|$. Likewise, $t_{0j}^k=t\times |l_0^k-l_j|$ and 
$t_{iT}^k=t\times |l_i-l_T^k|$  are the travelling times from the initial 
position to the bay of task $j$ and from the bay of task $i$ to the final 
position, respectively. If $l_T^k=0$ then the latter travelling time is set to 
$0$. Cranes are numbered according to their initial position so that the leftmost crane in the vessel is crane $1$ 
and the rightmost one is crane $q$.

In scheduling terminology, and disregarding cranes interferences due to crossing 
or safety margins, the problem corresponds to a minimum makespan scheduling 
problem with parallel identical machines (cranes) and precedence constraints, 
which is known to be \textbf{NP}-Hard in the strong sense, provided that more 
than two machines, non-preemption or non-uniform processing times are given 
\cite{pinedo2008}.

Next, we introduce the MIP model we use for the QCSP. To this end, consider the 
following set of variables:
\begin{itemize}
\item $x_{ij}^k\in \{0,1\}~\forall i\in \Omega^{0},~j\in \Omega^{T},~k\in K$. 
$x^k_{ij}=1$ if and only if task $j$ follows $i$ in the task sequence of quay 
crane $k$, otherwise $x^k_{ij}=0$. If $i=0$, $j$ is the first task performed by 
crane $k$ and, if $j=T$, $i$ is the last task performed by crane $k$;
\item $y_{ik}\in \{0,1\}~\forall i\in \Omega,~k\in K$. $y_{ik}=1$ if and only if task $i$ is assigned to quay crane $k$;
\item $z_{ij}\in \{0,1\}~\forall i,j\in \Omega$. $z_{ij}=1$ if and only if task j starts later than the completion of task $i$ (i.e. $i$ is completed before $j$ starts);
\item $D_i$ is the completion time of task $i\in \Omega$;
\item $C_k$ is the completion time of quay crane $k\in K$;
\item $W$ is the makespan, the completion time of the vessel.
\end{itemize}

For the sake of simplicity, the following notation will be used in the mathematical context throughout the text, where $S\subset \Omega \cup \{0,T\}$:
\begin{itemize}
\item $x^k(S)=\sum_{i,j \in S}x^k_{ij}$
\item $x^k(i,S)=\sum_{S \in j}x^k_{ij}$
\item $x^k(S,i)=\sum_{S \in j}x^k_{ji}$
%\item $\Omega^0=\Omega \cup \{0\}$
%\item $\Omega^T=\Omega \cup \{T\}$
\end{itemize}

The MIP model for the QCSP proposed by \cite{kim2004} with the developments and 
modifications by \cite{moccia2006} is the following:
\begin{align}
&\text{min } \alpha_1 W + \alpha_2 \sum_{k\in K}C_k \label{mdl:fo}\\
x^k(0,\Omega^T) = 1  &\quad \forall k \in K \label{mdl:dgr0}\\
x^k(\Omega^0,T) = 1  &\quad \forall k \in K \label{mdl:dgrT}\\
y_{ik} = x^k(i,\Omega^T) &\quad \forall i \in \Omega, \forall k \in K \label{mdl:out1}\\ 
y_{ik} = x^k(\Omega^0,i) &\quad \forall i \in \Omega, \forall k \in K \label{mdl:out2}\\
\sum_{k \in K} y_{ik} = 1 &\quad \forall i \in \Omega \label{mdl:flow}\\
D_i+t_{ij}+p_j-D_j \leq M(1-x_{ij}^k) &\quad \forall i,j \in \Omega, \forall k \in K \label{mdl:sch1}\\
D_i+p_j-D_j \leq M(1-z_{ij}) &\quad \forall i,j \in \Omega, l_i\neq l_j \label{mdl:sch2}\\
D_j-p_j-D_i \leq Mz_{ij} &\quad \forall i,j \in \Omega, l_i\neq l_j \label{mdl:sch3}\\
\sum_{v=1}^k y_{jv} + \sum_{v=k}^q y_{iv} \leq 1+z_{ij}+z_{ji} &\quad \forall i,j \in \Omega, l_i< l_j,\forall k\in K \label{mdl:cross}\\
D_i+p_j-D_j+\sum_{k\in K}\sum_{u\in \Omega^0, l_u \neq l_i} tx^k_{uj} \leq M(1-z_{ij}) &\quad \forall i,j \in \Omega, l_i=l_j \label{mdl:cross1}\\
D_j-p_j-D_i-\sum_{k\in K}\sum_{u\in \Omega^0, l_u \neq l_i} tx^k_{uj} \leq Mz_{ij} &\quad \forall i,j \in \Omega, l_i=l_j \label{mdl:cross2}\\
z_{ij}+z_{ji} = 1 &\quad \forall (i,j) \in \Psi \setminus \Phi \label{mdl:sch4}\\
z_{ij}=1, z_{ji}=0 &\quad \forall (i,j) \in \Phi \label{mdl:sch5}\\
r_k-D_j+t_{0j}^k+p_j \leq M(1-x_{0j}^k) &\quad \forall j\in \Omega,~\forall k\in K \label{mdl:sch6}\\
D_j+t_{jT}^k-C^k \leq M(1-x_{jT}^k) &\quad \forall j \in \Omega,~\forall k \in K \label{mdl:sch7}\\
C^k \leq W &\quad \forall k \in K \label{mdl:sch8}\\
x_{ij}^k, y_{ik}, z_{ij} \in \{0,1\} & \quad \forall i,j\in \Omega,~\forall k \in K
\end{align}

The objective is to minimize a weighted sum of the makespan and the sum of 
completion times for each crane. In related works using the instance set 
proposed by \cite{kim2004}, minimizing the makespan is the primary objective 
($\alpha_1 \gg \alpha_2)$. Constraints \eqref{mdl:dgr0}-\eqref{mdl:flow} are the 
routing constraints of the cranes. As noted by \cite{moccia2006}, it is possible 
for a crane $k$ to leave its initial state, $0$, and go directly to the final 
state, $T$, without performing any task. Constraints \eqref{mdl:cross} are the 
non-crossing constraints proposed by \cite{moccia2006}. They were designed to prevent 
crane crossing forcing that in case tasks $i$ and $j$, $l_i<l_j$, are performed 
simultaneously, then the crane assigned to task $i$ must be lower than the crane 
assigned to task $j$. Inequalities \eqref{mdl:sch1} determine the completion 
times of tasks and eliminate subtours in the cranes routes. Inequalities 
\eqref{mdl:sch2}-\eqref{mdl:cross2} define variables $z_{ij}$ concerning the 
order in which tasks are performed. If tasks $i$ and $j$ cannot be processed 
simultaneously, i.e. $(i,j)\in \Psi \setminus \Phi$, then we must either have 
$i\prec j$ or $j \prec i$ by inequalities \eqref{mdl:sch4}. Constraints 
\eqref{mdl:sch5} define the precedence relationships among tasks. The start time 
of the first task performed by each crane is defined by inequalities 
\eqref{mdl:sch6}. The completion time of each crane is defined by inequalities 
\eqref{mdl:sch7} and the makespan of the schedule is determined by inequalities 
\eqref{mdl:sch8}. The constant M is a large constant, but in \cite{moccia2006} 
the authors propose several ways to reduce this value in each inequality.

\subsection{Crane interference}
As noted earlier, \cite{bier2009} observed that the previous model still lacks a 
correct treatment of crane interference constraints. In their revised model 
they introduced a suitable temporal distance between any two tasks involved in 
a problem. Let $\Delta_{ij}^{vw}$ denote the minimum time to elapse between the 
processing of tasks $i$ and $j$ if assigned to cranes $v$ and $w$, respectively. 
The correction proposed by \cite{bier2009} to address the crane interferences 
results from the replacement of constraints \eqref{mdl:cross}, 
\eqref{mdl:cross1} and \eqref{mdl:cross2} by \eqref{mdl:crossR}, 
\eqref{mdl:cross1R} and \eqref{mdl:cross2R}, respectively.
\begin{align}
y_{iv} + y_{jw} \leq 1+z_{ij}+z_{ji} &\quad (i,j,v,w) \in \Theta \label{mdl:crossR}\\
D_i+\Delta_{ij}^{vw}+p_j-D_j \leq M(3-z_{ij}-y_{iv}-y_{jw})&\quad (i,j,v,w) \in \Theta \label{mdl:cross1R}\\
D_j+\Delta_{ij}^{vw}+p_i-D_i \leq M(3-z_{ij}-y_{iv}-y_{jw})&\quad (i,j,v,w) \in \Theta \label{mdl:cross2R}
\end{align}

where $\Theta=\{(i,j,v,w)\in \Omega^2 \times K^2 | i<j, \Delta_{ij}^{vw}>0\}$ is the set of all two tasks and crane assignments that potentially lead to a crossing situation. If task $i$ is assigned to crane $v$ and task $j$ is assigned to crane $w$, then the left side of constraints \eqref{mdl:crossR} is two, hence tasks $i$ and $j$ are not allowed to be performed simultaneously, since $z_{ij}+z_{ji}=1$. If task $j$ starts after the completion of task $i$, i.e. $z_{ij}=1$, then constraints \eqref{mdl:cross1R} insert a suitable temporal distance $\Delta_{ij}^{vw}$ between the completion of $i$ ($D_i)$ and before the start of $j$ ($D_j-p_j$) so that the assigned cranes can move without interference. The case $z_{ji}=1$ is analogous and handled by constraints \eqref{mdl:cross2R}.

\subsection{Crane limits}
The allocation of cranes to ships must abide to several constraints such as 
technical data about cranes and ships and the accessibility of cranes to a 
berth. Considering the integration of quayside problems, the assignment must 
reflect vessels' neighbour berths and all ships moored in the terminal. The 
temporal distance included in the model ensures a sufficient time span to elapse 
between the processing of tasks $i$ and $j$, allowing a safe movement of the 
cranes. On the other hand, a large safety margin might lead to a situation 
in which cranes can be positioned outside the limits of the ships' bays. This 
situation may not be desirable, for example, if the crane interferes with 
another crane operating in an adjacent ship.

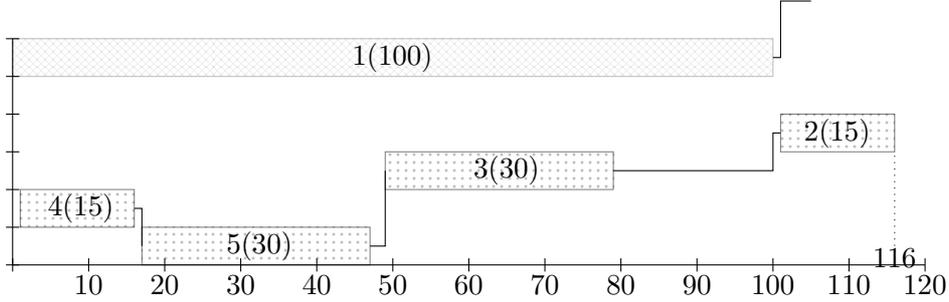
\begin{figure}
\centering
\caption{Example of a crane outside the limits of the vessel bays}
\label{fig:ex1}
\begin{tikzpicture}[]
    % Axis
    \draw (0,0) -- (12,0);
    \draw (0,0) -- (0,3);
    % Note that the snaked line is drawn to 11.1 to force
    % TikZ to draw the final tick.
    \draw[snake=ticks,segment length=0.5cm] (0,0) -- (0,3);
    \draw[snake=ticks,segment length=1cm] (0,0) -- (12,0);
    \node[below] at (1,0) {10};
    \node[below] at (2,0) {20};
    \node[below] at (3,0) {30};
    \node[below] at (4,0) {40};
    \node[below] at (5,0) {50};
    \node[below] at (6,0) {60};
    \node[below] at (7,0) {70};
    \node[below] at (8,0) {80};
    \node[below] at (9,0) {90};
    \node[below] at (10,0) {100};
    \node[below] at (11,0) {110};
    \node[below] at (12,0) {120};  

    \filldraw[pattern=crosshatch, opacity=0.25] (0,2.5) rectangle (10,3); %task 1
    \filldraw[pattern=dots, opacity=0.5] (0.1,0.5) rectangle (1.6,1); %task 4
    \filldraw[pattern=dots, opacity=0.5] (1.7,0) rectangle (4.7,0.5); %task 5
    \filldraw[pattern=dots, opacity=0.5] (4.9,1) rectangle (7.9,1.5); %task 3
    \filldraw[pattern=dots, opacity=0.5] (10.1,1.5) rectangle (11.6,2); %task 2
    
    \draw[] (10,2.75)-- (10.1,2.75)-- (10.1,3.5)-- (10.5,3.5);
    \draw[] (1.6,0.75)-- (1.7,0.75)-- (1.7,0.25);
    \draw[] (4.7,0.25)-- (4.9,0.25)-- (4.9,1.25);
    \draw[] (7.9,1.25)-- (10,1.25)-- (10,1.75)-- (10.1,1.75);
    
    \node[] at (5,2.75) {1(100)};
    \node[] at (0.9,0.75) {4(15)};
    \node[] at (3.25,0.25) {5(30)};
    \node[] at (6.5,1.25) {3(30)};
    \node[] at (10.85,1.75) {2(15)};
    
    \draw[dotted] (11.6,0) -- (11.6,1.5);
    \node[] at (11.6, 0.1) {116};
\end{tikzpicture}
\end{figure}

One way to ensure that cranes will operate within the limits of the vessel bays 
is to impose that some bays can only be visited by certain cranes. For instance, 
if two cranes are assigned to a ship then the rightmost crane can not reach the 
first $1+\delta$ bays, otherwise the leftmost crane would need to be positioned 
outside the vessel. Figure \ref{fig:ex1} illustrates such a situation, the travel 
time is one time unit. Tasks 1 and 2 cannot be processed simultaneously due the 
safety margin $\delta=2$. If task 2 is assigned to crane 2, the safety margin 
will only be respected if crane 1 is moved to the left of bay 1 while task 2 is 
processed. The optimal makespan is 116 in this case. If the limits are imposed, 
task 2 can only be processed by crane 1, the optimal schedule 
changes, and the optimal makespan is increased by one unit of time.

The leftmost and the rightmost bays in which crane $k\in K$ can operate are 
defined by $l_m^k=(k-1)(\delta+1)+1$ and $l_M^k=B-(K-k)(1+\delta)$, 
respectively. To impose 
such operational limits for the cranes in the model, the following constraints 
are included:
\begin{align}
\sum_{\substack{i\in \Omega \\ l_i<l_m^k \\ l_i>l_M^k}} y_{ik} =0 &\quad \forall k \in K \label{mdl:limits}
\end{align}

In constraints \eqref{mdl:limits}, if $l_i$ is not within the range of crane $k$, then task $i$ cannot be assigned to that crane, i.e. $y_{ik}=0$.

\section{Decomposition algorithm}
\label{sec:alg}
The QCSP can be decomposed into three consecutive steps. First, a crane must be 
assigned to each task. At this stage, 
%the compatibility \review{(note: change name)} task-crane 
constraints \eqref{mdl:limits} must be satisfied to avoid pushing cranes outside the ship boundaries. Then, the tasks assigned to each crane must be sequenced. Each sequence must obey the task precedence constraints. The sequences produced in this stage determine the route (the sequence of bays) each crane would follow to complete its assigned tasks in case no interference with other crane exists. Finally, a starting time is set to each task considering the route established in the previous step, the precedence constraints for tasks assigned to different cranes, the interference constraints and the safety margin. This last step also computes the solution makespan.

In this work, we propose an iterative decomposition approach tackling the first 
and second steps as a master sub-problem and the third step as a slave 
sub-problem. The master sub-problem, which we call the routing problem, determines the route of each crane and the slave sub-problem, which we call the 
scheduling problem, determines the completion time for each task given the 
routes for each crane.

The routing stage consists in solving a Vehicle Routing Problem (VRP) 
considering distinct initial and final depots and precedence constraints. The 
objective function we consider is the minimization of the longest route, 
considering the processing times of the tasks and the time to move from one bay 
to another. We propose a basic branch-and-cut algorithm for this sub-problem. 
Observe that the solution cost of this problem is a lower bound for the makespan 
of the full problem. 

Once crane routes are known, the scheduling problem can be formulated as a 
polynomially sized integer linear programming problem and be solved by any 
integer programming solver via branch and bound. This sub-problem is always 
feasible 
% \review{(given a sufficiently large $M$?)} 
and its solution cost is an upper bound for the makespan of the full 
problem.

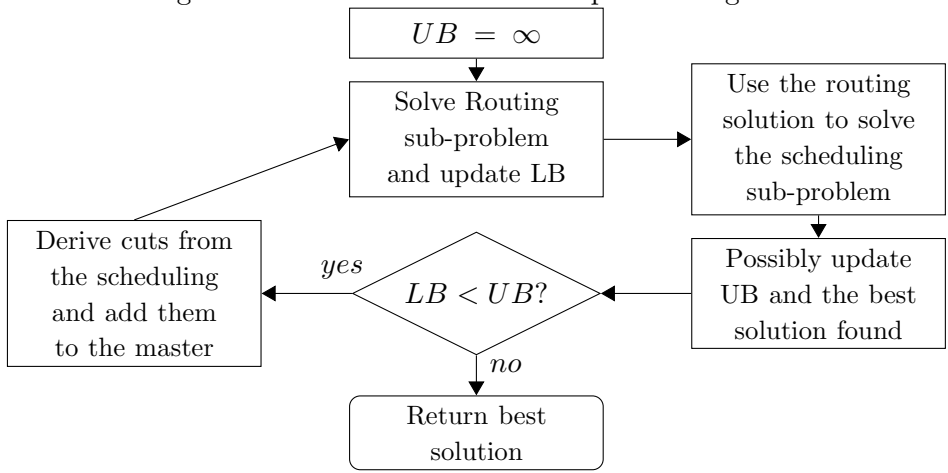
\begin{figure}
\caption{Flowchart for the decomposition algorithm}
\label{fig:flow}
% Start the picture
\begin{tikzpicture}
[
    >=triangle 60,              % Nice arrows; your taste may be different
    start chain=going below,    % General flow is top-to-bottom
    node distance=3mm and 45mm, % Global setup of box spacing
    every join/.style={norm},   % Default linetype for connecting boxes
]
% ------------------------------------------------- 
% A few box styles 
% <on chain> *and* <on grid> reduce the need for manual relative
% positioning of nodes
\tikzset{
  base/.style={draw, on chain, on grid, align=center, minimum height=4ex},
  proc/.style={base, rectangle, text width=8em},
  test/.style={base, diamond, aspect=2, text width=5em},
  term/.style={proc, rounded corners},
  % coord node style is used for placing corners of connecting lines
  coord/.style={coordinate, on chain, on grid, node distance=6mm and 25mm},
  % nmark node style is used for coordinate debugging marks
  nmark/.style={draw, cyan, circle, font={\sffamily\bfseries}},
  % -------------------------------------------------
  % Connector line styles for different parts of the diagram
  norm/.style={->, draw, black},
  free/.style={->, draw, lcfree},
  cong/.style={->, draw, lccong},
  it/.style={font={\small\itshape}}
}
% -------------------------------------------------
% Start by placing the nodes
\node [] (r) at (0,-5) {};
\node [proc] (p0) {$UB=\infty$};
\node [proc, join] (p1) {\small{Solve Routing sub-problem and update LB}};
%\node [proc, join, right= of p1] (p2) {Fix $z_{ij}$ for tasks $i,j$ assigned to the same crane};
\node [proc, join, right= of p1] {\small{Use the routing solution to solve the scheduling sub-problem}};
\node [proc, join] (p3)          {\small{Possibly update UB and the best solution found}};
\node [test, join, left= of p3] (t1) {$LB<UB$?};
\node [proc,  left= of t1] (p4) {\small{Derive cuts from the scheduling and add them to the master}};
\node [term, below= of r] (f1)  {\small{Return best solution}};

\path (p4.north) to node [near start, xshift=1em] {} (p0);
  \draw [->] (p4.north) -- (p1.west);
\path (t1.west) to node [near start, xshift=-1em] {$yes$} (p1);
  \draw [->] (t1.west) -- (p4);  
\path (t1.south) to node [near start, xshift=1em] {$no$} (f1);
  \draw [->] (t1.south) -- (f1);

\end{tikzpicture}
\end{figure}

In our approach the master problem is responsible for providing routes for the 
slave problem and for updating the problem lower bound. The slave problem 
updates the problem upper bound (cost of best known solution) and modifies the 
master problem in order to avoid the generation of the same or similar routes
that cannot further improve the best know solution in future iterations. When 
the lower bound computed by the master problem reaches the cost of the best 
known solution, the algorithm stops returning the current best known solution 
that is proved to be optimal at this point. The flowchart in Figure 
\ref{fig:flow} illustrates how the proposed method works.

\subsection{Combinatorial Benders' Cuts}
Benders' decomposition method \cite{benders1962} is a classical approach for dealing with large scale optimization problems. Basically, the method consists in splitting the MIP formulation into a master and a slave sub-problem. A solution to the master sub-problem provides a lower bound to the full problem and a solution to the slave sub-problem provides an upper bound and either optimality or feasibility cuts which are added to the master. The full problem is solved iteratively until the gap between the lower and upper bounds is sufficient small. Generally, the structure of the problem leads to a partitioning in which the master sub-problem works in the space of the complicating variables (for instance, integer ones) and the slave is a linear programming problem.

More recently, \cite{codato2006} developed a Benders' like method aimed to remove the model dependency on big-M values usually introduced in MIP formulations. In the combinatorial Benders' method, the linking between the integer and the continuous variables is assumed to be due to a set of constraints involving only one integer variable multiplied by the big-M coefficient. The master sub-problem is an integer program and the cuts returned by the slave sub-problem are purely combinatorial inequalities separated through the solution of a irreducible infeasible subsystem (IIS) of the slave constraints. For example, if $\textbf{x}\in [0,1]^n$ is a vector of binary variables in the master and a subset of these variables, say $C$, induces an IIS of the slave sub-problem given a solution $\bar{\textbf{x}}$, then not all variables 
%indexed by 
in $C$ can assume their actual value. In this case, the following  combinatorial Benders' cut is added to the master sub-problem:
\begin{align}
\sum_{i\in C:\bar{\textbf{x}}_i=0}\textbf{x}_i +\sum_{i\in C:\bar{\textbf{x}}_i=1}(1-\textbf{x}_i)\geq& 1 \label{comb}
\end{align}

% Applying the method proposed in \cite{codato2006} for the QCSP model presented 
% in Section \ref{sec:math} would require some reformulations. Note that 
% constraints \eqref{mdl:crossR}, \eqref{mdl:cross1R} and \eqref{mdl:cross2R} do 
% not conform with the structure of the liking inequalities. 
% Moreover, if binary variables $z_{ij}$ are defined in the master, inequalities 
% \eqref{mdl:crossR} are not sufficient to determine the order in which tasks $i$ 
% and $j$ are processed when assigned to different cranes. Note that in these 
% inequalities if task $i$ is assigned to crane $v$ and $j$ is assigned to crane 
% $w$, then the left side of constraints \eqref{mdl:crossR} is two, hence tasks 
% $i$ and $j$ are not allowed to be performed simultaneously, since we must have 
% $z_{ij}+z_{ji}=1$. However, the routing problem does not provide information 
% whether $z_{ij}=1$ or $z_{ji}=1$.\review{(improve this...)}

In our approach, we leave the $z_{ij}$ variables in the scheduling sub-problem 
and, therefore, we have a MIP instead of an LP sub-problem as is the case in the 
combinatorial Benders's approach. Nevertheless, if $i$ and $j$ are assigned to 
the same crane, the tasks sequence obtained after solving the routing problem 
allows to fix the value of $z_{ij}$. The binary variables not fixed in the 
scheduling problem are those $z_{ij}$ for which tasks $i$ and $j$ are assigned 
to different cranes. In our computational experiments, this MIP sub-problem 
could be solved with little computational effort.

\subsection{The master sub-problem}
The objective in the MIP formulation for the QCSP is a weighted sum of the 
makespan, $W$, and crane completion times, $C_k$. In this work, as in other 
works dealing with the QCSP, it is assumed that $\alpha_1 \gg \alpha_2$. In 
particular, we set $\alpha_2=0$. Therefore, the objective is to minimize the 
makespan. In the master problem, we define a cost for each routing variable 
$x_{ij}^k$ so that the solution of the routing problem provides a lower bound 
for the makespan. Let the cost $c_{ij}^k$ be such that:
\begin{align*}
c_{0i}^k = r_k + t_{0i}^k + p_i & \quad i \in \Omega,~k \in K\\
c_{ij}^k = t_{ij} + p_j & \quad i,j \in \Omega,~k \in K\\
c_{iT}^k = t_{iT}^k & \quad i \in \Omega,~k \in K
\end{align*}
that is, $c^k_{0i}$ is the cost (the travelling and processing time) of crane 
$k$ leaving its initial position, going to the bay where task $i$ is located and 
processing it. The cost $c_{ij}^k$ is the travelling time of crane $k$ moving from 
the bay $l_i$ to the bay $l_j$ and the processing time of task $j$ (note that this cost does not depend on the crane $k$). Finally, the cost $c_{iT}^k$ is the travelling time from bay $l_i$ to the final bay $l_T^k$ 
of crane $k$.

The cost vector $\textbf{c}$ defines a completion time for each crane and the objective function in the routing sub-problem is to minimize the maximum crane completion time:
\begin{align}
\text{min}~\eta&\label{mdl:foR}\\
 \sum_{j \in \Omega} c_{0j}^kx_{0j}^k +
 \sum_{i,j \in \Omega} c_{ij}^kx_{ij}^k +
 \sum_{i \in \Omega} c_{iT}^kx_{iT}^k & \leq \eta & \forall k \in K \label{mdl:minmax}
\end{align}

The routing sub-problem consists of \eqref{mdl:foR}, \eqref{mdl:minmax}, \eqref{mdl:dgr0}-\eqref{mdl:flow} and \eqref{mdl:limits}. Note that, if no interactions occur among cranes (i.e. there are no idle times due to crane interferences or crossings) then the optimal value of $\eta$ must be equal to the makespan obtained in the scheduling sub-problem.

\subsubsection{Solving the master sub-problem}
We resort to a branch-and-cut scheme to solve the routing sub-problem to 
optimality. Before start, we add to the master problem an initial pool of valid 
inequalities comprising of all sub-tour elimination constraints (SEC) 
$x^k(S)\leq |S| -1$ for $|S| = 2$. Also, we include in the pool the 
\textit{precedence cycle breaking inequalities}, proposed by \cite{balas1995} in 
the context of the Precedence Constrained Asymmetric Traveling Salesman 
Problem. These inequalities are added for each pair of precedences 
$(i_1,j_1),~(i_2,j_2)\in \Phi$ so the inequality becomes 
$x^k_{i_1j_2}+x^k_{j_1i_2}\leq 1~\forall k \in K$.

Since inequalities \eqref{mdl:sch1} are not present in the master, sub-tours might occur in the routing solution. When all variables $x_{ij}^k$ are integer, we separate violated SEC by identifying the connected components in the supporting graph of the solution. Observe that if an integer solution contains any sub-tour, then the support graph of the solution is induced by $q$ paths from $0$ to $T$ and by cycles $C_1,...,C_l\subset \Omega$ where each $C_i,~1\leq i\leq l$ defines a violated sub-tour elimination inequality.

Since 
%variables $z_{ij}$ are not defined in the master sub-problem and 
the order in which tasks are executed are determined in the routing 
problem, precedence relations need to be dealt with in this stage using only 
the arc variables $x^k_{ij}$. We distinguish between two cases, namely, when the two 
tasks in a precedence relationship are performed by the same crane and when they 
are performed by different cranes.

Precedence violations in tasks performed by the same crane are identified by 
adapting a SEC lifting proposed by \cite{balas1995} in the context of the 
Precedence Constrained Traveling Salesman Problem. If $i\prec j$, for a subset 
$S$ such that $0,j\in S$ and $i\notin S$, the sub-tour elimination constraints 
can be lifted in the following way: $x^k(S)\leq |S|-2$. For the QCSP 
formulation, that is not true if tasks $i$ and $j$ are processed by different 
cranes, say, $k_i$ and $k_j$. Crane $k_j$ can arrive at bay $l_j$ before crane 
$k_i$ arrives at bay $l_i$, but need to wait until task $i$ is processed to 
continue and handle task $j$, thus $i\prec j$. We consider the following 
inequality:
\begin{align}
x^k(S)+y_{ik}\leq |S|-1&\quad \forall S,~0,j\in S~i\notin S,~(i,j)\in \Phi,~k\in K \label{mdl:prec}
\end{align}
When $i$ and $j$ are processed by the same crane $k$, $y_{ik}=1$ and the inequality becomes $x^k(S)\leq |S|-2$, which is violated if the crane processes task $j$ before task $i$. If that is not the case, the inequality is trivially satisfied, since $S$ must not contain a sub-tour.

\begin{figure}
\caption{Infeasible routing solution due to precedence relations.}
\label{fig:precS}
\centering
\begin{tikzpicture}
[place/.style={circle,draw=black!80,fill=black!60,thick,
inner sep=0pt,minimum size=3mm},
depot/.style={rectangle,draw=black!80,fill=black!60,thick,
inner sep=0pt,minimum size=3mm}]
\node at ( 3,1.5) [depot, label=left:{$0$}] (v1) {};
\node at ( 3,0.5) [depot, label=left:{$0$}] (v2) {};
\node at ( 5,0.5) [place, label=above:{$j$}] (v3) {} edge[<-] (v2);
\node at ( 6,1.5) [place, label=above:{$m$}] (v4) {}  edge[<-] (v1);
\node at ( 7,0.5) [place, label=above:{$l$}] (v6) {} edge[<-] (v3);
\node at ( 9,0.5) [depot, label=right:{$T$}] (v5) {} edge[<-] (v6);
\node at ( 8,1.5) [place, label=above:{$i$}] (v7) {}  edge[<-] (v4);
\node at ( 9,1.5) [depot, label=right:{$T$}] (v8) {} edge[<-] (v7);

\node at (1,1.5) [] {Crane 1};
\node at (1,0.5) [] {Crane 2};
\end{tikzpicture}
\end{figure}
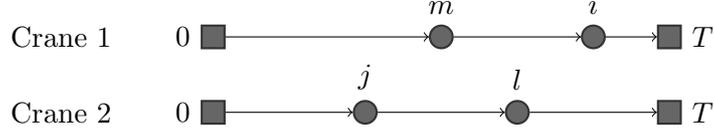

When no precedence constraints are violated within each route, it is still 
possible that the routing obtained leads to an infeasible solution. In Figure 
\ref{fig:precS} we depict an example of a infeasible solution, where $q=2$ and 
$(i,j),(l,m)\in \Phi$. Tasks $i$ and $m$ are handled by crane $1$ and tasks $j$ 
and $l$ by crane $2$. If $m$ is processed before $i$ by crane $1$ and task $j$ 
is processed before $l$ by crane $2$, then this solution is infeasible, since 
at least one of the precedences $(i,j)$ or $(l,m)$ is violated.

This situation can be generalized for more than two cranes. For instance, if 
$q=3$ and $(i,j),~(l,m),~(p,q)\in \Phi$ then a solution in which task $m$ is 
processed before $i$ by a crane $k_1$, task $q$ is processed before $l$ by crane 
$k_2$ and task $j$ is processed before $p$ by crane $k_3$ would lead to the 
same kind of precedence infeasibility. If an integer solution does not violate 
any inequality \eqref{mdl:prec}, then we try to find $g \le q$ ordered pairs 
$(i_1,j_1),...,(i_g,j_g)\in \Phi$, where tasks $j_m$ and $i_{(m+1 \mod g)+1}$ 
are processed by crane $k_m$, leading to the aforementioned situation (note 
that all tasks need to be different). 
The following inequality is then added:
\begin{align}
x^{k_1}(j_1,S_1)+x^{k_1}(S_1)+...+x^{k_g}(j_g,S_g)+x^{k_g}(S_g)\leq 
|S_1|+...+|S_g|-1 \label{mdl:precS}
\end{align}
where $S_m$ consists of a path from $j_m$ (excluded) to $i_{(m+1 \mod g)+1}$ 
(included). That is, inequality \eqref{mdl:precS} impose that at least one of 
the paths from $j_i$ to $i_{(j)}$ need to be changed to circumvent the precedence 
violations.

% Instead of inequalities \eqref{mdl:precS}, another two inequalities proposed in \cite{balas1995}, namely, the predecessor and the successor inequalities could be used. Since precedence relations are defined for tasks in the same bay, in most of the routings obtained, the two tasks in a precedence relation are processed by the same crane and violations like these seldom occur. We opted for inequalities \ref{mdl:precS} since they are easier to define and to separate.

\subsection{The slave sub-problem}
With the task sequence for each crane obtained after solving the routing sub-problem, each variable $z_{ij}$ is fixed for tasks $i$ and $j$ assigned to the same crane, i.e. if $i$ is processed before $j$ by a given crane, then $z_{ij}=1$ and $z_{ji}=0$. The routing solution $\bar{\textbf{x}}$ and the partial fixing $\bar{\textbf{z}}$ are used to build the scheduling sub-problem, consisting of the objective function \eqref{mdl:fo}, where $\alpha_2=0$, and \eqref{mdl:sch1}--\eqref{mdl:sch3}, \eqref{mdl:sch4}--\eqref{mdl:sch8}, \eqref{mdl:crossR}--\eqref{mdl:cross2R}.

We solve the scheduling sub-problem for the given routing solution and the
corresponding optimal makespan $W^*$ is obtained. If the solution makespan is smaller 
than the best solution found so far, then we update the best solution. In 
this case, the 
combinatorial cut sent to the master is a simple 'no good' 
cut derived from \eqref{comb} where $C$ is the subset of those $x_{ij}^k$ 
variables which are equal to $1$ in the solution, that is, the cut
\begin{align}
\sum_{x \in C}x \leq |C|-1 \label{mdl:nogood}
\end{align}
eliminates the current routing from the solution space of the subsequent routing sub-problems.

\subsubsection{Cut generation}
The no good inequality can be strengthened to cut more than one single routing 
from the solution space. Suppose that tasks $i$ and $j$ are in the same bay and 
that they are processed sequentially, say, $i$ and then $j$, by the same crane in the 
routing solution. A solution in which $j$ is processed immediately before 
$i$ by the same crane would lead to a schedule with exactly the same makespan (provided that there is no precedence relation between tasks $i$ and $j$). 
If $p$ is the task processed immediately before $i$ and $q$ the task processed 
immediately after $j$, that is $(p,i),(i,j),(j,q) \in C$, then sum 
$x^k_{pj}+x^k_{ji}+x^k_{iq}$ can be added to the left side of 
\eqref{mdl:nogood}. More generally, if $S$ is a subset of tasks in the same bay 
processed by a crane $k$, $p$ is the task processed immediately before tasks in 
$S$ and $q$ is the task processed immediately after, then the variables representing the path of crane $k$ from $p$ to $q$ in $C$ can be replaced by $x^k(p,S)+x^k(S)+x^k(S,q)$ in the left side of \eqref{mdl:nogood}. 
Note that the right-hand side of \eqref{mdl:nogood} remains the same, that is, 
the sum of all crane's routing size subtracted by one.

%In the combinatorial Benders' method, violated cuts are generated solving an IIS for the given slave sub-problem. Observe that, in the scheduling formulation, any routing solution is feasible if the value $M$ is sufficiently large. On the other hand, suppose that the optimal makespan is $W^*$ after solving the first scheduling sub-problem. This value could be used to replace $M$ in the following iterations ($W^*+\Delta_{ij}^{vw}$ in inequalities \eqref{mdl:cross1R} and \eqref{mdl:cross2R}) and infeasible scheduling sub-problems could be derived. One of the sources of infeasibility is due to the fact that the precedence relations for tasks processed by different cranes might lead to an optimal makespan which could be larger than the best makespan obtained in previous iterations. Instead of solving an IIS, we identify such cases using a sufficiently large value $M$ such that the scheduling sub-problem is always feasible.

\begin{figure}
\caption{Routing solution for which tasks $i$ and $j$ ($i\prec j$) are processed by different cranes.}
\label{fig:cut}
\centering
\begin{tikzpicture}
[place/.style={circle,draw=black!80,fill=black!60,thick,
inner sep=0pt,minimum size=3mm},
depot/.style={rectangle,draw=black!80,fill=black!60,thick,
inner sep=0pt,minimum size=3mm}]
\node at ( 3,1.5) [depot, label=left:{$0$}] (v1) {};
\node at ( 3,0.5) [depot, label=left:{$0$}] (v2) {};
\node at ( 5,0.5) [place, label=above:{$j$}] (v3) {} edge[<-, dashed] (v2);
\node at ( 6,1.5) [place, label=above:{$i$}] (v4) {}  edge[<-] (v1);
\node at ( 9,0.5) [depot, label=right:{$T$}] (v6) {} edge[<-] (v3);
\node at ( 9,1.5) [depot, label=right:{$T$}] (v7) {} edge[<-, dashed] (v4);

\node at (1,1.5) [] {Crane 1};
\node at (1,0.5) [] {Crane 2};
\end{tikzpicture}
\end{figure}
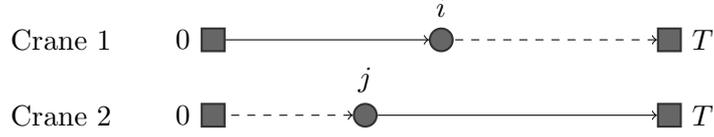

Consider now two tasks $i$ and $j$, $(i,j)\in \Phi$, which are processed by different cranes, say, $k_i$ and $k_j$, respectively, as shown in Figure \ref{fig:cut}. Since $i\prec j$, tasks in the path from $j$ to $T$ cannot be processed before task $i$ is handled. Thus, if crane $k_i$ takes at least $s_i$ units of time from $0$ to $i$ (including travel and processing times, excluding cranes interferences) and crane $k_j$ takes at least $s_j$ units of time from $j$ to $T$, this solution is sub-optimal since $UB^*<s_i+s_j$, where $UB^*$ is the best makespan found so far. In that case, rather than add to the master sub-problem a no good inequality \eqref{mdl:nogood}, the cut can be improved considering in $C$ only $x_{ij}^k$ variables in the paths of cranes $k_i$, from $0$ to $i$, and $k_j$, from $j$ to $T$ (the solid lines in Figure \ref{fig:cut}).

The cut can be further strengthened to eliminate more sub-optimal routes. Let $S_i$ be the set of tasks in the path from $0$ to $i$ of crane $k_i$, such that $0\in S_i$, $i\notin S_i$ and let $S_j$ be the set of tasks in the path from $j$ to $T$ of crane $k_j$, $j\notin S_j$, $T \in S_j$. If the sum of the processing times of tasks in $S_i$, $S_j$, $i$ and $j$ is greater than $UB^*$, then no matter in which order the tasks in $S_i$ and $S_j$ are processed, the route is still sub-optimal and the following cut can be included:
\begin{align}
x^{k_i}(S_i)+x^{k_i}(S_i,i)+x^{k_j}(j,S_j)+x^{k_j}(S_j)\leq |S_i|+|S_j|-1\quad (i,j)\in \Phi \label{mdl:Sset}
\end{align}

Inequality \eqref{mdl:Sset} cuts from the routing solution space all routes (i.e. task-to-crane assignments and sequencing) containing a subpath from $0$ to $i$ through $S_i$, for crane $k_i$, and a subpath from $j$ to $T$ through $S_j$, for crane $k_j$. For each pair of tasks $(i,j)$ such that $i\prec j$, a cut can be obtained by a simple inspection of the routing solution. If no such cut could be found, then we resort to a no good inequality.

\section{Computational study}
\label{sec:resul}
Our algorithm was implemented in C++ and uses the B\&C framework provided by Gurobi 5.63. As we consider instances with integer-valued costs, we set the absolute MIP gap tolerance and the absolute objective difference cut-off parameters to 0.9999. The tests were run on an Intel Core i5 CPU 760 @ 2.80GHz with 8GB RAM memory.

The first benchmark data set used in our experiments is a suite of instances introduced by \cite{kim2004}. The instances are numbered from 13 to 49 and are divided into four sets, namely, $A$ (instances $13$ to $22$), $B$ ($23$ to $32$), $C$ ($33$ to $41$) and $D$ ($42$ to $49$). Table \ref{tab:inst} shows the description of each instance set. The ready times for each crane, $r_k$, are set to $0$, the travelling time between two adjacent bays, $t$, is set to one unit of time and the safety margin distance, $\delta$, is set to one empty bay between any two cranes. The problem size varies from 10 to 25 tasks and from two to three quay cranes assigned to the vessel. The subset comprising the 37 instances numbered from 13 to 49 are tackled in \cite{kim2004}, \cite{moccia2006} and \cite{samarra2007}. For 28 of these instances, the optimal solutions are known. The revised model proposed by \cite{bier2009} corrected the optimal value of instance 22 and the best know solution for instance 42.
\begin{table}
\caption{Description of the instance set}
\label{tab:inst}
\centering
\begin{tabular}{ccccc}
\hline
Instance set & A & B & C & D \\
\hline
\hline
Tasks  & 10 & 15 & 20 & 25\\
Cranes & 2  & 2  & 3  & 3\\
\hline
\end{tabular}
\end{table}

\subsection{Limits on Crane Movements}
First, we evaluate the impact of imposing limits on the bays that each crane can reach, through equations \eqref{mdl:limits}, on the instance set. Let Formulation \textit{F1} be the full model, composed by \eqref{mdl:fo}--\eqref{mdl:sch3}, \eqref{mdl:sch4}--\eqref{mdl:sch8} and \eqref{mdl:crossR}--\eqref{mdl:limits}, and let \textit{F2} be the same model without equations \eqref{mdl:limits}, that is, the corrected formulation for the QCSP proposed by \cite{bier2009}.
\begin{table}
\caption{Comparison between formulations \textit{F1} and \textit{F2}.}
\label{tab:form12}
\centering
\begin{tabular}{cccrccrc}
\hline
\multicolumn{2}{c}{Instance} & & \multicolumn{2}{c}{F1} & & \multicolumn{2}{c}{F2}\\
\cmidrule{1-2} \cmidrule{4-5} \cmidrule{7-8}
Name & Makespan  & & Time(s) & Gap(\%) &  & Time(s) & Gap(\%)\\
\hline
$k13$ & 151 & & 1.45 & 0.00 & & 15.89 & 0.00\\   
$k14$ & 182 & & 1.28 & 0.00 & & 5.70 & 0.00\\	 
$k15$ & 171 & & 0.17 & 0.00 & & 2.33 & 0.00\\
$k16$ & 104 & & 0.14 & 0.00 & & 7.29 &	0.00\\
$k17$ & 151 & & 0.67 & 0.00 & & 7.12 & 0.00\\
$k18$ & 125 & & 0.65 & 0.00 & & 9.53 & 0.00\\
$k19$ & 181 & & 1.90 & 0.00 & & 14.301 & 0.00\\
$k20$ & 133 & & 0.73 & 0.00 & & 10.224 & 0.00\\
$k21$ & 155 & & 0.23 & 0.00 & & 0.481 & 0.00\\
$k22$ & 180 & & 6.86 & 0.00 & & 21.252 & 0.00\\
\hline
$k23$ & 192 & & 7.02 & 0.00 & & 45.165 & 0.00\\
$k24$ & 222 & & 21.34 & 0.00 & & 63.231 & 0.00\\
$k25$ & 246 & & 36.05 & 0.00 & & 63.716 & 0.00\\
$k26$ & 213 & & 7.58 & 0.00 & & 20.176 & 0.00\\
$k27$ & 219 & & 2.12 & 0.00 & & 21.531 & 0.00\\
$k28$ & 177 & & 12.99 & 0.00 & & 41.329 & 0.00\\
$k29$ & 269 & & 12.63 & 0.00 & & 86.427 & 0.00\\
$k30$ & 297 & & 2.45 & 0.00 & & 11.928 & 0.00\\
$k31$ & 190 & & 278.36 & 0.00 & & 1264.79 & 0.00\\
$k32$ & 197 & & 17.83 & 0.00 & & 53.88 & 0.00\\
\hline
$k33$ & 201 & & 308.29 & 0.00 & & 2634.83 & 0.00\\
$k34$ & \textbf{239} & & \textbf{347.94} & \textbf{0.00} & & 7200 & 0.00\\
$k35$ & \textbf{228} & & \textbf{437.82} & \textbf{0.00} & & 7200 & 0.01\\
$k36$ & \textbf{226} & & \textbf{4301.32} & \textbf{0.00} & & 7200 & 0.01\\
$k37$ & 170 & & 543.91 & 0.00 & & 6683.52 & 0.00\\
$k38$ & n.a. & & 7200 & 0.01 & & 7200 & 0.02\\
$k39$ & n.a. & & 7200 & 0.01 & & 7200 & 0.02\\
$k40$ & \textbf{188} & & \textbf{3563.96} & \textbf{0.00} & & 7200 & 0.03\\
$k41$ & \textbf{196} & & \textbf{2035.43} & \textbf{0.00} & & 7200 & 0.01\\
\hline
$k42$ & n.a & & 7200 & 0.03 & & 7200 & 0.03\\
$k43$ & n.a & & 7200 & 0.02 & & 7200 & 0.03\\
$k44$ & n.a & & 7200 & 0.01 & & 7200 & 0.02\\
$k45$ & n.a. & & 7200 & 0.02 & & 7200 & 0.03\\
$k46$ & n.a. & & 7200 & 0.01 & & 7200 & 0.01\\
$k47$ & \textbf{264} & & \textbf{479.12} & \textbf{0.00} & & 7200 & 0.01\\
$k48$ & n.a. & & 7200 & 0.02 & & 7200 & 0.03\\
$k49$ & n.a. & & 7200 & 0.02 & & 7200 & 0.03\\
\hline
\textit{Average} &  & & \textbf{2397.70} &  & & 3218.51 & \\
\hline
\end{tabular}
\end{table}

For each instance, we allow both formulations a maximum computation time of two hours. We report the results in Table \ref{tab:form12}. For each instance, we show the optimal makespan obtained with both formulations (the makespan values obtained with each model could be different, but they where all equal for the considered instance set) or n.a. (not achieved), if neither formulation was capable of finding an optimal solution for the given instance. The highlighted rows correspond to instances for which only formulation \textit{F1} was able to provide an optimal solution within the 2-hours limit. Observe that for all instances, the execution times with \textit{F1} are smaller than with formulation \textit{F2}. For this instance set, maintaining the cranes within the limits of the vessel did not increase the makespan values obtained with formulation \textit{F1}. Moreover, eliminating those variables corresponding to the movement of a crane to a bay which will lead to violation of the limits proved to be an effective approach to solve the formulation.

\subsection{Results with the Decomposition Approach}
In order to evaluate the effectiveness of the decomposition approach introduced in Section \ref{sec:alg}, we compare the results obtained with formulation \textit{F1} against the results obtained using the proposed approach on the hardest instance sets, $C$ and $D$. The maximum computation time is set to two hours. In Table \ref{tab:decomp}, we depict the results obtained with formulation \textit{F1} and with the decomposition approach for each instance in the set. In the column Makespan, we report the optimal makespan value obtained for each instance, or n.a. (not achieved) if neither algorithm was able to yield an optimal makespan for the given instance within the time limit. Column Time shows the computational time spent for each method and column Gap gives the final optimality gap of formulation \textit{F1}. Column lb/ub shows the lower and upper bounds obtained with the decomposition approach.
\begin{table}
\caption{Comparison of results obtained with the complete formulation and with the decomposition algorithm.}
\label{tab:decomp}
\centering
\begin{tabular}{cccrccrc}
\hline
\multicolumn{2}{c}{Instance} & & \multicolumn{2}{c}{Formulation \textit{F1}} & & \multicolumn{2}{c}{Decompostion}\\
\cmidrule{1-2} \cmidrule{4-5} \cmidrule{7-8}
Name & Makespan  & & Time(s) & Gap(\%) &  & Time(s) & lb/ub\\
\hline
$k33$ & 201 & & 308.29 & 0.00 & & 5.75 & 201/201\\
$k34$ & 239 & & 347.94 & 0.00 & & 11.06 & 239/239\\
$k35$ & 228 & & 437.82 & 0.00 & & 6.46 & 228/228\\
$k36$ & 226 & & 4301.32 & 0.00 & & 479.74 & 226/226\\
$k37$ & 170 & & 543.91 & 0.00 & & 45.53 & 170/170\\
$k38$ & n.a. & & 7200 & 0.01 & & 7200 & $205$/$218$\\
$k39$ & \textbf{171} & & 7200 & 0.01 & & 1531.15 & 171/171\\
$k40$ & 188 & & 3563.96 & 0.00 & & 7200 & $187$/$350$\\
$k41$ & 196 & & 2035.43 & 0.00 & & 4080.23 & 196/196\\
\hline
$k42$ & n.a & & 7200 & 0.03 & & 7200 & $189$/$191$\\
$k43$ & n.a & & 7200 & 0.02 & & 7200 & $290$/$\inf$\\
$k44$ & n.a & & 7200 & 0.01 & & 7200 & $273$/$367$\\
$k45$ & \textbf{278} & & 7200 & 0.02 & & 363.79 & 278/278\\
$k46$ & \textbf{230} & & 7200 & 0.01 & & 249.47 & 230/230\\
$k47$ & 264 & & 479.12 & 0.00 & & 11.28 & 264/264\\
$k48$ & n.a. & & 7200 & 0.02 & & 7200 & $211$/$349$\\
$k49$ & n.a. & & 7200 & 0.02 & & 7200 & $296$/$503$\\
\hline
\textit{Average} &  & & 4655.67& -- & & 3592.05 & --\\
\hline
\end{tabular}
\end{table}

Observe \ref{tab:decomp} that the decomposition approach reduced the computational times for the sets containing the larger instances. Two instances in set $D$, $k45$ and $k46$, could only be solved by the decomposition approach. In set $C$, instance $k39$ is only solved by our decomposition method. On the other hand, instance $k40$ is only solved with formulation \textit{F1}. For the remaining instances of the set, computational times are decreased one order of magnitude, and only for instance $k41$ our decomposition method is slower than solving it with \textit{F1}.

An optimal solution for some instances in sets $C$ and $D$ can be obtained by the decomposition approach after solving the first master sub-problem. This situation occurs for instances $k33$, $k34$, $k35$, $k46$ and $k47$, for which the decomposition method solves just one master sub-problem and it is three orders of magnitude faster when compared with the formulation \textit{F1}. Whereas formulation \textit{F1} could not solve $k46$, the decomposition approach converges after just one routing problem. In these cases, the master solution obtained is an unidirectional schedule without any crossing situations among cranes and, therefore, no idle times or temporal distances are added after solving the scheduling sub-problem to allow safe crane movements. The optimal makespan returned by the scheduling sub-problem is the same as the optimal cost of the master sub-problem, the lower and upper bounds are equal after the first iteration of the decomposition and the algorithm finishes with an optimal schedule.

Even when more than one routing sub-problem need to be solved, as is the case for instances $k36$, $k37$, $k39$, $k41$ and $k45$, computational times are substantially decreased. The decomposition method could solve $k45$ after adding seven cuts to the master sub-problem (one 'no good' cut \eqref{mdl:nogood} and six cuts \eqref{mdl:Sset}), but formulation \textit{F1} did not converge within the time limit on this instance.

Finally, we observe that solving the master sub-problem is the most time-consuming task in the decomposition. For instance $k43$, the time limit is exceeded during the solution of the first sub-problem. Even though the scheduling sub-problem includes binary variables, fixing variables $z_{ij}$ when tasks $i$ and $j$ are processed by the same crane allows the scheduling sub-problem to be solved with little computational effort.

\subsection{Literature results}
\subsubsection{Kim and Park instances}
Next, we compare the results obtained with our decomposition approach with another exact method from the literature, the B\&C algorithm of \cite{moccia2006}, and with the results obtained with the heuristic procedure of \cite{bier2009} using the previous set of benchmark instances. We note that the B\&C method by \cite{moccia2006} did not address crane interference correctly and that the model proposed by \cite{bier2009} corrected the optimal makespan value for instance $k22$ and the best known solution for instance $k42$. 

The corrected makespan for instance $k22$ (180) obtained with the revised QCSP model from \cite{bier2009} was proved to be optimal in this work by both solving formulation \textit{F1} and applying the decomposition approach.

Our approach is the first exact method to prove the optimal makespan for three instances. For the set $C$, our decomposition method could solve two instances, $k39$ and $k41$, not previously solved to optimality by the B\&C in \cite{moccia2006}. Just for one instance previously solved ($k40$) our decomposition algorithm did not converge within the two-hour limit. For set $D$, our decomposition approach was the first to prove the optimal solution for instance $k45$ and the other two instances previously solved in the literature ($k46$ and $k47$) were also solved to optimality. In these three cases, the optimal makespan returned by our decomposition method is the same as the best solution found by the heuristic in \cite{bier2009}. In fact, the optimal schedules returned by our approach for these three instances are unidirectional schedules. 

All instances in set $A$ and $B$ are solved to optimality by both the B\&C of \cite{moccia2006} and a solution of the same quality is yield by the heuristic proposed by \cite{bier2009}. Solving these instances using Formulation \textit{F1}, we also could prove the optimal makespan value for each of them.

\subsubsection{Bierwirth and Meisel instances}
%\review{
In the benchmark instances of \cite{kim2004}, the number of tasks always equals the number of bays and, consequently, the average number of tasks per bay (\textit{task-per-bay ratio}) is 1. Moreover, instances with more than 25 tasks are unrealistic, considering that the largest vessels nowadays contain around 24 bays \cite{chen2014}. Also, as pointed by \cite{bier2011}, at most four containers groups are assigned to a bay, leading to a unique sequential processing order. Finally, the assignment of tasks to bays and the processing times of tasks results in a workload which is not set into relation with the capacity of the bays. Considering these facts, \cite{bier2011} developed a new generation scheme to produce QCSP instances under a range of different parameters values, such as tasks-per-bay ratio, work-load distribution, precedence relation density, capacities of the bays, handling rates, among others.

We consider the instance set \textit{A} from the benchmark suite provided by \cite{bier2011}. The 70 instances consist of vessels with 10 bays, each with capacity of 200 containers and served by two cranes. The number of container groups (tasks) varies from 10 to 40.

In \cite{bier2011}, the 70 instances in set \textit{A} were handed to the unidirectional search heuristic (UDS) of \cite{bier2009} and to CPLEX 11 (with a maximum computational time of two hours) in order to assess the quality of the solutions obtained with the UDS. In our experiments, we also set a time limit of two hours for both solving the formulation \textit{F1} and the decomposition approach. In Table \ref{tab:instBM}, we report the results obtained and compare them with the results reported in \cite{bier2011} (Ref. Table 6). In bold, we highlight those  instances for which the model in \cite{bier2011} could not be solved within the time limit by CPLEX but were solved by formulation \textit{F1} and/or the decomposition approach. If the time limit is reached in both for a given instance, then we report $lb/ub$ in the correspondent line, where $lb$ is the best lower bound obtained by the decomposition approach and $ub$ is the best solution found by Gurobi within two hours.
%}

\begin{table}
\caption{Instance set A proposed by \cite{bier2011}.}
\label{tab:instBM}
\centering
\begin{tabular}{cccccccc}
\hline
No. & $n=10$ & $n=15$ & $n=20$ & $n=25$ & $n=30$ & $n=35$ & $n=40$ \\
\hline
1  & 520$^2$ & $\textbf{513}^{\diamond,1}$  & \textbf{508}$^2$ & \textbf{508}$^1$ & \textbf{506}$^1$ & 506$^1$ & \textbf{506/506}\\
2  & 508$^2$ & 507$^2$  & \textbf{509}$^2$ & 507$^2$ & 507/508 & 507/528 & 506$^1$\\
3  & 513$^2$ & \textbf{513/513} & \textbf{509}$^2$ & \textbf{507}$^2$ & 507/510 & 506$^1$ & 505$^1$\\
4  & 510$^2$ & $509/513^\ddagger$ & \textbf{509}$^1$ & 507/508$\ddagger$ & \textbf{507}$^1$ & 506/510 & --/507\\
5  & $\textbf{514}^{\diamond,1}$ & 507$^2$ & 506$^2$ & \textbf{507}$^1$ & 506$^1$ & --/510 & 506$^1$\\
6  & 513$^2$ & 508$^2$ & \textbf{508}$^1$ & \textbf{507}$^1$ & \textbf{506}$^1$ & 509/519 & \textbf{507/507}\\
7  & 511$^1$ & 507$^2$ & \textbf{507}$^2$ & 507/508 & 507/510 & 506/507 & $507/513\ddagger$\\
8  & 513$^2$ & 508$^2$ & \textbf{510}$^2$ & 507$^2$ & 506/510 & \textbf{506}$^1$ & 506$^1$\\
9  & 512$^1$ & 507$^2$ & \textbf{508}$^2$ & 506$^1$ & 506$^2$ & \textbf{506}$^1$ & 506$^1$\\
10 & 549$^1$ & 513$^1$ & 507$^1$ & 506$^2$ & 506$^1$ & 507/514 & \textbf{507/507}\\
\hline
\end{tabular}
    \begin{tablenotes}
      \small
      \item $\diamond:$ optimal schedule is not unidirectional.
      \item $\ddagger:$ previously solved to optimality in \cite{bier2011}.
      \item $1:$ \textit{F1} solved the instance first
      \item $2:$ Decomposition solved the instance faster than decomposition
    \end{tablenotes}
\end{table}

%\review{
Overall, the formulation \textit{F1} and the decomposition approach were able to provide optimality certificates for 23 open instances. Note that, for four cases included in those 23, the best $lb$ returned by the decomposition approach equals the makespan of the best solution ($ub$) obtained with formulation \textit{F1}, hence this solution is indeed an optimal schedule. For two instances, namely, instance $5$ for $n=10$ and instance $1$ for $n=15$, the solutions obtained have vessel handling times (VHT) $514$ and $513$, respectively, one unit time less than the VHT obtained with UDS (and equal the lower bound obtained by CPLEX). Indeed, for these two instances, the schedules obtained with both the formulation \textit{F1} and the decomposition are not unidirectional. In fact, for most of the instances, the optimal schedule obtained is unidirectional whenever the crane with completion time determining the VHT is unidirectional. For the two aforementioned instances this is not the case, an indication that a unidirectional schedule is not optimal.

In Table \ref{tab:instBM}, a superscript $2$ on an entry indicates that the computational time required by the decomposition approach is smaller when compared to solve \textit{F1} for a given instance, otherwise we use $1$. Computational time required to solve instances up to $n=25$ tasks could be improved with the decomposition approach when compared to solve formulation \textit{F1} An optimal solution could be found after solving a few master problems, the most cumbersome task in the algorithm, in some cases just one iteration was sufficient. For the larger instances, solving the master sub-problem proved to be very hard. Recall that, in this new proposed set of instances, since the \textit{task-per-bay} is greater than $1$, more precedence relations need to be satisfied when compared to the instances in the benchmark of \cite{kim2004}. For two instances, the time limit was reached solving the first master sub-problem (a $-$ on Table \ref{tab:instBM} indicates this situation). For this reason the sets of larger instances proposed in \cite{bier2011}, that were previously only tackled by heuristics, were not considered in this work. Finally, two instances solved to optimality in \cite{bier2011} could not be solved by either \textit{F1} nor the decomposition approach.
%}

\section{Conclusions and future research}
\label{sec:concl}
The quay crane scheduling problem is an important problem faced by container terminals in ports for the transshipment of cargo between container vessels and other modes of transportation. A more efficient use of the cranes available for loading and unloading the containers can decrease the time vessels spend in port and improve the throughput of the port.

In this work, we improved a MIP model for the QCSP in the literature by restricting the cranes to move outside the boundaries of the vessel. The solution obtained with this new model might better reflect the real operations of vessel in the port. Moreover, by eliminating variables that cannot be part of a solution with this new restriction, the model could be solved with less computational effort.

We proposed an algorithmic approach for solving the MIP formulation decomposing it into a master vehicle
routing problem and a corresponding slave scheduling problem. Our algorithm was capable of providing optimality certificates for three unsolved instances in a classical benchmark. A more recent suite was also used and 23 open instances could be solved to optimality. In particular, we could provide non-unidirectional schedules with a smaller makespan than the previous best unidirectional schedule known for two instances.

Future work includes the development of stronger cuts that can be derived from the scheduling obtained after solving the schedule sub-problem. We also aim to tackle a model which integrates all the quay side operations of a vessel, including the berth allocation, quay crane assignment and the quay crane scheduling, using similar decomposition techniques.

\section*{Acknowledgements}
The work presented here was carried out within the Port-Ship Coordinated Planning project, supported by the Norwegian Research Council under project number 227084/O70.

%\bibliographystyle{spbasic}    % basic style, author-year citations
%\bibliographystyle{spmpsci}     % mathematics and physical sciences
%\bibliographystyle{spphys}     % APS-like style for physics
\bibliographystyle{plain}
\bibliography{references}       % name your BibTeX data base

\end{document}